\documentclass[11pt]{article}
\usepackage{etoolbox}
\newtoggle{arxiv}
\toggletrue{arxiv}
\usepackage{makeidx}  % allows for indexgeneration
\usepackage{url}
\usepackage{graphicx}

\usepackage{amsfonts,amsmath}
\def\bE{{\mathbb E}}
\def\bR{{\mathbb R}}
\def\bS{{\mathbb S}}
\def\bZ{{\mathbb Z}}
\def\cC{{\mathcal C}}
\def\cF{{\mathcal F}}
\def\cH{{\mathcal H}}
\def\cL{{\mathcal L}}
\def\cO{{\mathcal O}}
\def\gC{{\mathfrak C}}

\def\diff{{\mathrm d}}
\def\targetSet{\Theta}
\def\sourceSet{\Upsilon}
\def\cot{\hat}
\def\<{\langle} \def\>{\rangle}
\def\sm{\setminus}
\def\ve{\varepsilon}

\renewcommand\subset{\subseteq}
\begin{document}
%
%\frontmatter          % for the preliminaries
%
%\pagestyle{headings}  % switches on printing of running heads
%\addtocmark{Geodesic methods with constraints} % additional mark in the TOC
%

%\mainmatter              % start of the contributions
%
%\title{Computation of optimal flag capture trajectories, surveillance optimization}
%\subtitle{Automatic differentiation of non-holonomic fast marching.}
\title{Automatic differentiation of non-holonomic fast marching for computing most threatening trajectories under sensors surveillance}
%
%\titlerunning{Hamiltonian Mechanics}  % abbreviated title (for running head)
%                                     also used for the TOC unless
%                                     \toctitle is used
%
\iftoggle{arxiv}{
\author{Jean-Marie Mirebeau\footnote{University Paris-Sud, CNRS, University Paris-Saclay} \and Johann Dreo\footnote{THALES Research \& Technology}}

}{ % llncs
\author{Jean-Marie Mirebeau\inst{1} \and Johann Dreo\inst{2}}
\authorrunning{Jean-Marie Mirebeau et al.} % abbreviated author list (for running head)
%
%%%% list of authors for the TOC (use if author list has to be modified)
%\tocauthor{Ivar Ekeland, Roger Temam, Jeffrey Dean, David Grove,
%Craig Chambers, Kim B. Bruce, and Elisa Bertino}
%
\institute{University Paris-Sud, CNRS, University Paris-Saclay,\\
\email{jean-marie.mirebeau@math.u-psud.fr}
\and
THALES Research \& Technology\\
\email{johann.dreo@thalesgroup.com}
}
}

\maketitle              % typeset the title of the contribution

\begin{abstract}
We consider a two player game, where a first player has to install a surveillance system within an admissible region. The second player needs to  enter the the monitored area, visit a target region, and then leave the area, while minimizing his overall probability of detection. 
Both players know the target region, and the second player knows the surveillance installation details.

Optimal trajectories for the second player are computed using a recently developed variant of the fast marching algorithm, which takes into account curvature constraints modeling the second player vehicle maneuverability. 
The surveillance system optimization leverages a reverse-mode semi-automatic differentiation procedure, estimating the gradient of the value function related to the sensor location in time $\cO(N \ln N)$.
\iftoggle{arxiv}{
\paragraph{keywords: }Anisotropic fast marching, motion planning, sensor placement, optimization
}{
}
\end{abstract}

\section{Introduction}

This paper presents a proof of concept numerical implementation of a motion planning algorithm related to a two player game.
A first player selects, within a admissible class $\Xi$, an integral cost function on paths, which takes into account their position, orientation, and possibly curvature. 
The second player selects a path, within an admissible class $\Gamma$, with prescribed endpoints and an intermediate keypoint. 
The players objective is respectively to maximize and minimize the path cost
\begin{align}
\label{eq:ToyProblem}
	\gC(\Xi, \Gamma) &:= \sup_{\xi \in \Xi} \inf_{\gamma \in \Gamma } \gC(\xi,\gamma),&
	\text{where } \gC(\xi,\gamma) &:= \int_0^{T(\gamma)} \cC_\xi(\gamma(t), \gamma'(t), \gamma''(t)) \, \diff t,
\end{align}
where the path $\gamma$ is parametrized at unit euclidean speed, and the final time $T(\gamma)$ is free.
From a game theoretic point of view, this is a non-cooperative zero-sum game, where player $\Xi$ has no information and player $\Gamma$ has full information over the opponent's strategy.

The game \eqref{eq:ToyProblem} typically models a surveillance problem~\cite{Strode:multistatic}, and $\exp(-\gC(\Xi, \Gamma))$ is the probability for player $\Gamma$ to visit the intermediate keypoint without being detected by player $\Xi$ (by convention, players are identified with their sets of strategies). For instance player $\Gamma$ could be a robber trying to steal an object from a museum, and player $\Xi$ would be responsible for the placement of video surveillance cameras.
In another instance, player $\Xi$ is responsible for the installation of radar~\cite{BarbarescoMonnier:threatening} or sonar detection systems~\cite{Strode:multistatic}, and would like to prevent vehicles sent by player $\Gamma$ from spying on some objectives without being detected. 

The dependence of the cost $\cC_\xi$ w.r.t.\ the path tangent $\gamma'(t)$ models the variation of a measure of how detectable the target is (radar cross section, directivity index, etc.) w.r.t.\ the relative positions and orientations of the target and sensor.
The dependence of $\cC_\xi$ on the path curvature $\gamma''(t)$ models the airplane maneuverability constraints, such as the need to slow down in tight turns~\cite{Duits2017:ReedsShepp}, or even a hard bound on the path curvature~\cite{Dubins1957:Car}.

Strode~\cite{Strode:multistatic} have shown the interplay of motion planning and game theory in a similar setting, on a multistatic sonar network use case, but using isotropic graph-based path planning. The same year, Barbaresco~\cite{Barbaresco:threatAniso} used  fast-marching for computing threatening paths toward a single radar, but without taking into account curvature constraints and without considering a game setting. 
% have shown that anisotropic % can be useful
% It is not clear at all to me that Barbaresco uses anisotropic fast marching in this paper. 

The main contributions of this paper are as follows:
%\vspace{-\topsep}
\begin{enumerate}
	\item {\em Anisotropy and curvature penalization}:
	Strategy optimization for player $\Gamma$ is an optimal motion planning problem, with a known cost function. This is addressed by numerically solving a generalized eikonal PDE posed on a two or three dimensional domain, and which is strongly anisotropic in the presence of a curvature penalty and a detection measurement that depends on orientation.
	%  and sophisticated
	A Fast-Marching algorithm, relying on recent adaptive stencils constructions, based on tools from lattice geometry, is used for that purpose~\cite{Duits2017:ReedsShepp,Mirebeau2017:VR1,Mirebeau2017:Curv}.
	In contrast, the classical fast marching method~\cite{Rouy1992:Viscosity} used in~\cite{Benmansour2010:Derivatives} is limited to cost functions $\cC_\xi(\gamma(t))$ independent of the path orientation~$\gamma'(t)$ and curvature~$\gamma''(t)$.
	
	\item {\em Gradient computation for sensors placement}:
	Strategy optimization for player $\Xi$ is typically a non-convex problem, to which various strategies can be applied, yet gradient information w.r.t.\ the variable $\xi \in \Xi$ is usually of help.
	For that purpose, we implement efficient differentiation algorithms, forward and reverse, for estimating the gradient of the value function of player $\Xi$
	\begin{align}
	\label{eqdef:ObjGradient}
		&\nabla_\xi \gC(\xi,\Gamma), &
		\text{where } \gC(\xi, \Gamma) := \inf_{\gamma \in \Gamma} \gC(\xi, \gamma).
	\end{align}
	Reverse mode differentiation reduced the computation cost of $\nabla_\xi \gC(\xi, \Gamma)$ from $\cO(N^2)$, as used in \cite{Benmansour2010:Derivatives}, to $\cO(N\ln N)$, where $N$ denotes the number of ization points of the domain. As a result, we can reproduce examples from \cite{Benmansour2010:Derivatives} with computation times reduced by several orders of magnitude, and address complex three dimensional problems.
\end{enumerate}

Due to space constraints, this paper is focused on problem modeling, ization and numerical experiments, rather than on mathematical aspects of wellposedness and convergence analysis.
Free and open source codes for reproducing (part of) the presented numerical experiments are available on the first author's webpage\footnote{\url{github.com/Mirebeau/HamiltonianFastMarching}}.

\section{Mathematical background of trajectory optimization}
\label{sec:MathBackground}
We describe in this section the PDE formalism, based on generalized eikonal equations, used to compute the value function $\min_{\gamma \in \Gamma} \gC(\xi,\gamma)$ of the second player, where $\xi$ is known and fixed. Their discretization is discussed in \S \ref{sec:Discretization}. We distinguish two cases, depending on wether the path cost function $\cC_\xi(x,\dot x,\ddot x)$ appearing in \eqref{eq:ToyProblem} depends on the last entry $\ddot x$, i.e.\ on path curvature.

\subsection{Curvature independent cost}
Let $\Omega \subset \bE := \bR^2$ be a bounded domain, and let the source set $\sourceSet$ and target set $\targetSet$ be subsets of $\Omega$.
For each $x \in \Omega$, let $\Gamma_x$ denote the set of all paths $\gamma \in C^1([0,T], \Omega)$, where $T = T(\gamma)$ is free, such that 
\begin{align*}
	\gamma(0) &\in \sourceSet, &
	\gamma(T) &= x, &
	\text{ and } \forall t \in [0,T],\, \|\gamma'(t)\| &= 1.
\end{align*}
The problem description states that the first player needs to go from $\sourceSet$ to $\targetSet$ and back, hence its set of strategies is $\Gamma = \bigcup_{x \in \targetSet} \Gamma_x \times \Gamma_x$. Therefore 
\begin{align}
\label{eq:TargetCostSplit_NoCurv}
	\gC(\xi,\Gamma) &= 2\inf_{x \in \targetSet} u_\xi(x), &
	\text{where } u_\xi(x) &:= \inf_{\gamma \in \Gamma_x} \gC(\xi,\gamma).
\end{align}
Define the $1$-homogenous metric $\cF_\xi : \Omega \times \bE \to [0,\infty]$, the Lagrangian $\cL_\xi$ and the Hamiltonian $\cH_\xi$ by 
\begin{align*}
	\cF_\xi(x,\dot x) &:= \|\dot x\| \cC_\xi(x,\dot x/\|\dot x\|), &
	\cL_\xi &:= \frac 1 2 \cF_\xi^2, &
	\cH_\xi(x,\cot x) &:= \sup_{v\in \bE} \<\cot x,\dot x\> - \cL_\xi(x,\dot x).
\end{align*}
Under mild assumptions~\cite{Bardi1997:OptimalControl}, the function $u_\xi : \Omega \to \bR$ is the unique viscosity solution to a generalized eikonal equation
\begin{align*}
	\forall x \in \Omega\sm \sourceSet,\, \cH_\xi(x,\nabla_x u_\xi(x)) &= 1/2, &
	\forall x \in \sourceSet,\, u_\xi(x) = 0,
\end{align*}
with outflow boundary conditions on $\partial \Omega $. The discretization of this PDE is discussed in \S \ref{sec:Discretization}.
We limit in practice our attention to Isotropic costs $\cC_\xi(x)$, and Riemannian costs $\cC_\xi(x, \dot x) = \sqrt{\<\dot x,M_\xi(x) \dot x\>}$ where $M_\xi(x)$ is symmetric positive definite, for which efficient numerical strategies have been developed \cite{Mirebeau2014:FastMarching,Mirebeau2017:VR1}.

\subsection{Curvature dependent cost}

Let $\Omega \subset \bR^2 \times \bS^1$ be a bounded domain, within the three dimensional space of all positions and orientations. As before, let $\sourceSet,\targetSet \subset \Omega$. For all $x \in \Omega$ let $\Gamma^\pm_x$ be the collection of all $\gamma \in C^2([0,T], \Omega)$, such that $\eta := (\gamma,\pm\gamma')$ satisfies
\begin{align*}
	\eta(0) &\in \sourceSet, &
	\eta(T) &= x, &
	\text{ and } \forall t \in [0,T],\, \|\gamma'(t)\| &= 1.
\end{align*}
Here and below, the symbol ``$\pm$'' must be successively replaced with ``$+$'' and then ``$-$''.
Since the first player needs to go from $\sourceSet$ to $\targetSet$ and back, its set of strategies is $\Gamma = \bigcup_{x \in \targetSet} \Gamma_x^+ \times \Gamma_x^-$. 
Therefore 
\begin{align}
\label{eq:TargetCostSplit_Curv}
	\gC(\xi,\Gamma) &= \inf_{x \in \targetSet} u^+_\xi(x)+u_\xi^-(x), &
	\text{where } u_\xi^\pm(x) &:= \inf_{\gamma \in \Gamma_x} \gC^\pm(\xi,\gamma).
\end{align}
We denoted $\gC^\pm$ the path cost defined in terms of the local cost $\cC_\xi(p,\pm \dot p,\ddot p)$. 
Consider the $1$-homogeneous metric $\cF^\pm_\xi : {\mathrm T} \Omega \to [0,\infty]$,  defined on the tangent bundle to $\Omega \subset \bR^2 \times \bS^1$ by 
%Focusing on the case of $u^+_\xi$, we introduce the Lagrangian $\frac 1 2 \cL$ 
\begin{align*}
	\cF^\pm_\xi( (p,n), (\dot p,\dot n)) := 
	\begin{cases}
	+\infty  & \text{if } \dot p \neq \|\dot p\| n,\\
	\|\dot p\| \cC_\xi(p,\pm n,\dot n/\|\dot p\|) &\text{else},
	\end{cases}
\end{align*}
where $p \in \bR^2$, $n \in \bS^1$ is a unit vector, and the tangent vector satisfies $\dot p \in \bR^2$, $\dot n \perp n$.

Introducing the Lagrangian $\cL^\pm_\xi = \frac 1 2 (\cF^\pm_\xi)^2$ on ${\mathrm T} \Omega$, and its Legendre-Fenchel dual the Hamiltonian $\cH^\pm_\xi$, one can again under mild assumptions characterize $u_\xi^\pm$ as the unique viscosity solution to the generalized eikonal PDE $\cH^\pm_\xi(x,\nabla u^\pm_\xi(x)) = 1/2$ with appropriate boundary conditions~\cite{Bardi1997:OptimalControl}. In practice, we choose cost functions of the form $\cC_\xi(p,\dot p, \ddot p) = \cC^\circ_\xi(p, \dot p) \cC_*(|\ddot p|)$, where $\cC_*$ is the 
Reeds-Shepp car
%Euler elastica~\cite{Mumford1994:Elastica}, 
or Dubins car~\cite{Dubins1957:Car} curvature penalty, namely
\begin{align*}
	\cC_{\rm RS}(\kappa) &:= \sqrt{1+\rho^2\kappa^2}, &
%	\cC_{\rm RS}(\kappa) &:= 1+\rho^2\kappa^2, & 
	\cC_{\rm D}(\kappa) &:= 
	\begin{cases}
		1 &\text{if } |\rho \kappa|\leq 1,\\
		+\infty &\text{otherwise},
	\end{cases}
\end{align*}
where $\rho>0$ is a parameter which has the dimension of a curvature radius.
The Dubins car can only follows path which curvature radius is $\leq \rho$, whereas the Reeds-Shape car, in the sense of~\cite{Duits2017:ReedsShepp} and without reverse gear), can rotate into place if needed.
The Hamiltonian then has the explicit expression $\cH((p,n),(\cot p, \cot n)) = \frac 1 2 \cC^0_\xi(p,n)^{-2} \allowbreak \cH_*(n,(\cot p, \cot n))$ where
\begin{align*}
    \cH_{\rm RS} &= \frac 1 2 (\<\cot p,n\>_+^2+\|\cot n/\rho\|^2) &
%	\cH_{\rm E} &= \frac 1 8 (\<\cot p,n\>+\sqrt{\<\cot p,n\>^2+\|\cot n/\rho\|})^2, &
	\cH_{\rm D} &= \frac 1 2 \max \{0,\<\cot p, n\>+\|\cot n/\rho\|\}^2.
\end{align*}

\section{Discretization of generalized eikonal equations}
\label{sec:Discretization}
\begin{figure}
    \centering
    \def\h{2.5cm}
    \def\dir{Figures/Stencils/}
    \includegraphics[height=\h]{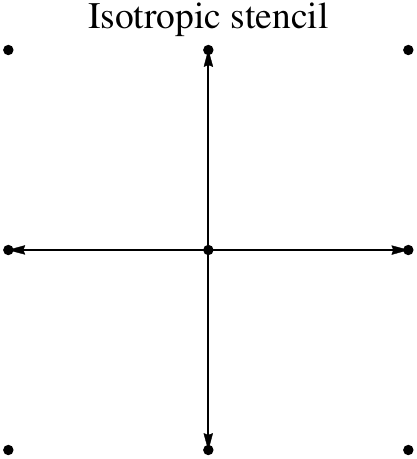}
    \includegraphics[height=\h]{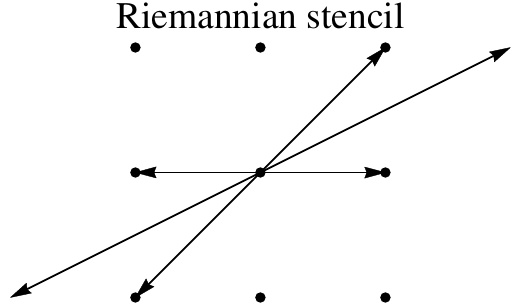}
    \includegraphics[height=\h]{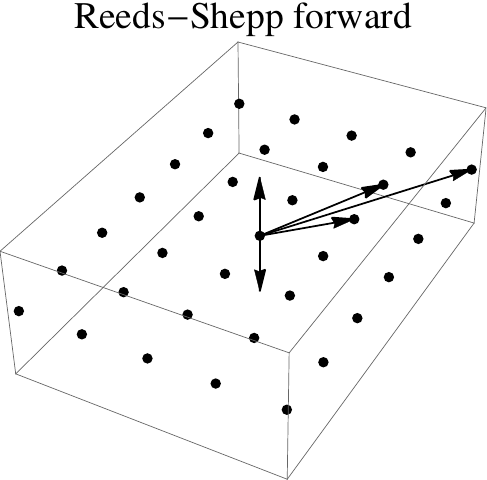}
    \includegraphics[height=\h]{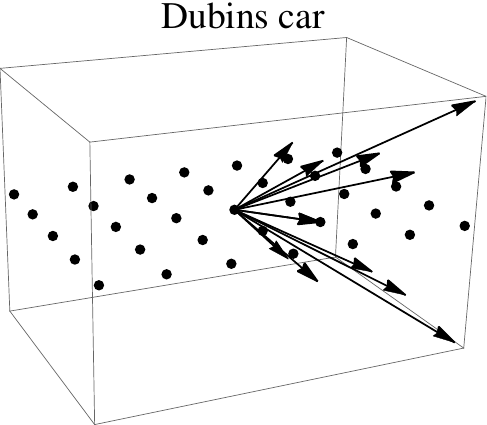}
    \caption{Examples of adaptive stencils used by the Fast-Marching solver}
    \label{fig:Stencils}
\end{figure}

We construct a discrete domain $X$ by intersecting the computational domain with an orthogonal grid of scale $h>0$
\begin{equation*}
	X = \Omega \cap (h \bZ)^d,
\end{equation*}
where $d=2$ for the curvature independent models, $d=3$ for the other models which are posed on $\bR^2 \times (\bR/2 \pi \bZ)$ ---using the angular parametrization $\bS^1 \cong \bR/2 \pi \bZ$ (in the latter periodic case, and $2\pi/h$ must be an integer). We design weights $c_\xi(x,y)$, $x,y \in X$ such that for any tangent vector $\dot x$ at $x$ one has
\begin{equation}
\label{eq:DiscreteHamiltonian}
	\cH_\xi(x,\dot x) \approx h^{-2} \sum_{y \in X} c^2_\xi(x,y)\<x-y, \dot x\>_+^2,
\end{equation}
where $a_+ := \max \{0,a\}$. (Expression \eqref{eq:DiscreteHamiltonian}, is typical although some models require a slight generalization.) The weights $c_\xi(x,y)$ are non-zero for only few $(x,y)\in X$ at distance $\|x-y\| = \cO(h)$, see Figure \ref{fig:Stencils}.
Their construction exploits the additive structure of the discretization grid $X$ and relies on techniques from lattice geometry~\cite{Schurmann2009:Quadratic}, see~\cite{Duits2017:ReedsShepp,Mirebeau2017:VR1,Mirebeau2017:Curv} for details.
The generalized eikonal PDE $\cH_\xi(x,\nabla_x u_\xi(x)) = 1/2$ is discretized as 
\begin{equation}
	\label{eq:Discretization}
	\sum_{y \in X} c^2_\xi(x,y) (U_\xi(x) - U_\xi(y))_+^2 = h^2/2
\end{equation}
with adequate boundary conditions. The solution $U_\ve : X \to \bR$ to this system of equations is computed in a single pass with $\cO(N \ln N)$ complexity using the Fast-Marching algorithm~\cite{Rouy1992:Viscosity}.

To be able to use the gradient to solve the problem~\eqref{eq:ToyProblem}, we need to differentiate the cost $\gC(\xi,\Gamma)$ w.r.t.\ the first player strategy $\xi \in \Xi$. In view of~\eqref{eq:TargetCostSplit_NoCurv} and~\eqref{eq:TargetCostSplit_Curv}, this only requires the sensitivity of the discrete solution values $U_\ve(x_*)$ at the few points points $x_* \in X \cap \targetSet$, w.r.t\ to variations in the weights $c_\ve(x,y)$, $x,y \in X$.
For that purpose we differentiate \eqref{eq:Discretization} w.r.t.\ $\xi$ and obtain
\begin{equation*}
	\sum_{y \in X} \omega_\xi(x,y) 
	\left(
	    \diff U_\xi(x)-\diff U_\xi(y) 
	  + (U_\xi(x) - U_\xi(y)) \diff \ln c_\xi(x,y)
	\right)
	= 0,
\end{equation*}
where $\omega_\xi(x,y) := c^2_\xi(x,y) (U_\xi(x) - U_\xi(y))_+$. This equation ties $\diff U_\ve(x)$ to those $\diff U_\xi(y)$, $y \in X$, such that $\omega_\xi(x,y)>0$, hence satisfying $U_\xi(x) < U_\xi(y)$ by construction. In the spirit of automatic differentiation by reverse accumulation~\cite{Griewank2008:Derivatives}, these relations are back-propagated from any $\targetSet$ keypoint $x_*$ to the $\sourceSet$ seeds, at a modest cost $\cO(N)$.

\section{Numerical results}
\label{sec:Num}

\begin{figure}
\def\w{4cm}
\includegraphics[width=\w]{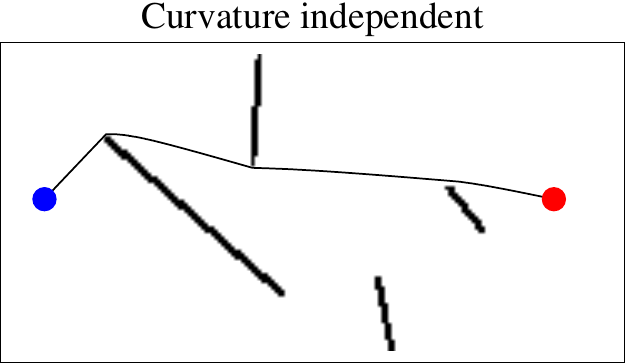}
\includegraphics[width=\w]{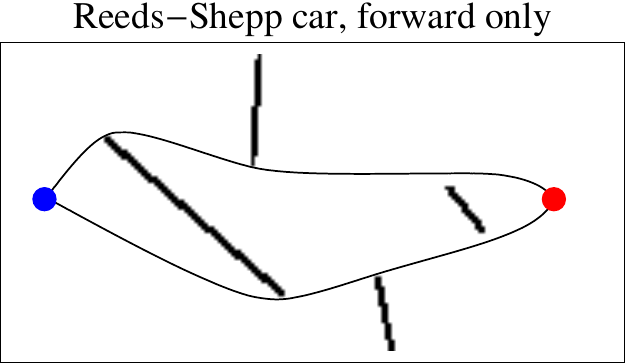}
\includegraphics[width=\w]{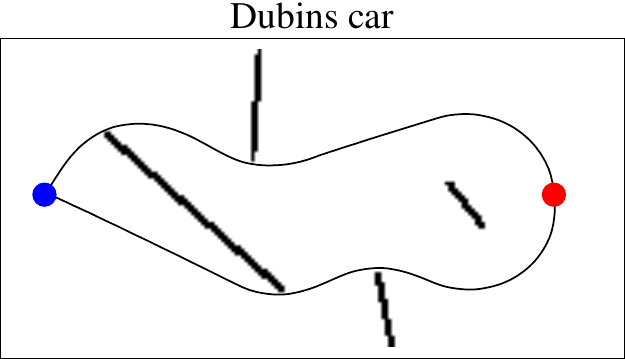}
\caption{No intervention from the first player. Shortest path from the blue point (left) to the red keypoint (right) and back.}
\label{fig:FreePaths}
\end{figure}

The chosen physical domain $R$ is the rectangle $[0,2]\times [0,1]$ minus some obstacles, as illustrated on Figure \ref{fig:FreePaths}. Source point is $(0.2,0.5)$ and  target keypoint $(1.8,0.5)$. 
The computational domain is thus $\Omega=R$ for curvature independent models and $\Omega = R\times \bS^1$ for curvature dependent models, which is discretized on a $180 \times 89$ or $180 \times 89 \times 60$ grid. 

\paragraph{No intervention from the first player.} The cost function is 
\begin{equation}
    \cC_\xi(p, \dot p, \ddot p) = \cC_*(|\ddot p|), 
\end{equation}
where $\cC_*(\kappa)$ is respectively $1$, $\sqrt{1+\rho^2\kappa^2}$ and $1+\infty\times \chi_{\rho\kappa>1}$, with $\rho := 0.3$.
The differences between the three models are apparent: the curvature independent model uses the same path forward and back; the Reeds-Shepp car spreads some curvature along the way but still makes an angle at the target point; the Dubins car maintains the radius of curvature below the bound $\rho$, and its trajectory is a succession of straight and circular segments.

\begin{figure}
\def\w{4cm}
\def\dir{Figures/Paint/Rectangle/}
\includegraphics[width=\w]{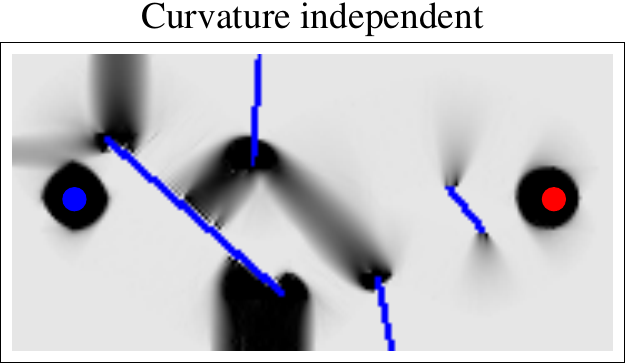}
\includegraphics[width=\w]{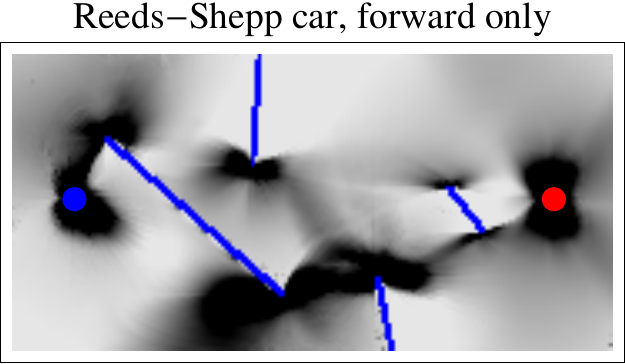}
\includegraphics[width=\w]{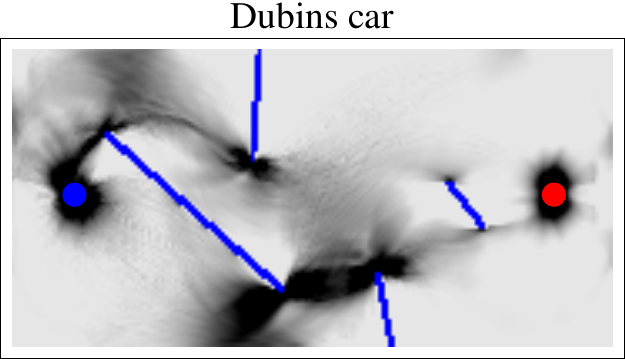}\\
\includegraphics[width=\w]{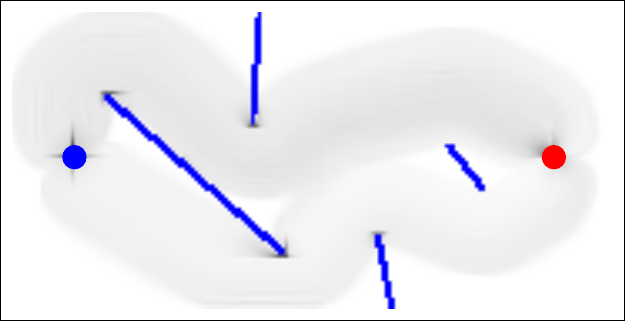}
\includegraphics[width=\w]{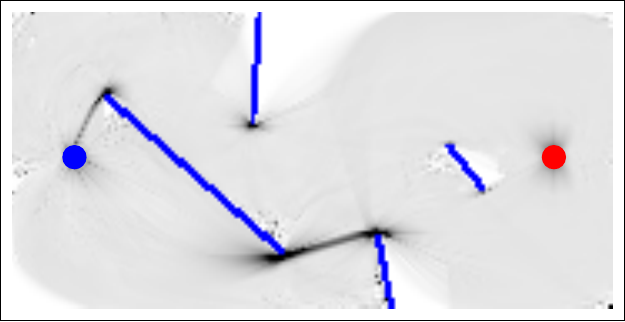}
\includegraphics[width=\w]{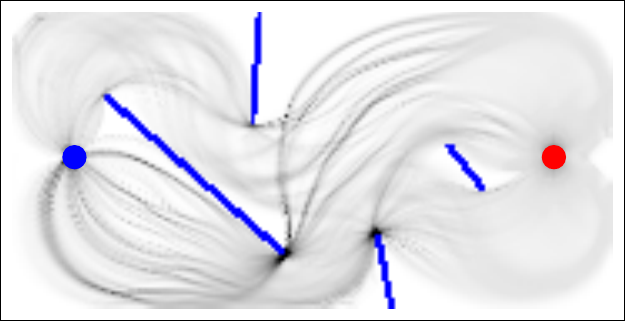}
\caption{Fresh paint based detection. Top: Optimal distribution of paint, to mark a path from the blue point (left) to the red keypoint (right) and back. Bottom: Geodesic density at the equilibrium.
}
\label{fig:Paint}
\end{figure}

\paragraph{Fresh paint based detection.} In this toy model, see Figure \ref{fig:Paint}, the first player spreads some fresh paint over the domain, and the second player is regarded as detected if he comes back covered in it from his visit to the keypoint. The cost function is 
\begin{equation}
    \cC_\xi(p,\dot p, \ddot p) = \xi(p) \cC_*(|\ddot p|)
\end{equation}
where $\xi : R \to \bR_+$ is the fresh paint density, decided by the first player, and $\cC_*(\kappa)$ is as above. For wellposedness, we impose a bounds on the density, namely $0.1 \leq \xi \leq 1$, and subtract the paint supply cost $\int_R \xi(p) \diff p$ to \eqref{eq:ToyProblem}.
The main interest of this specific game, also considered in \cite{Benmansour2010:Derivatives}, is that $\gC(\xi,\Gamma)$ is concave w.r.t.\ $\xi \in \Xi$. The optimal strategy for player $\Xi$ is: in the curvature independent case to make some ``fences'' of paint between close obstacles, and in the curvature penalized models to deposit paint at the edges of obstacles, as well as along specific circular arcs for the Dubins model. At the optimum for player $\Xi$, a continuum of paths have equal cost for player $\Gamma$, covering a large portion of the domain with a uniform density as illustrated on Figure \ref{fig:Paint} (bottom), see \cite{Benmansour2010:Derivatives} for details.

\begin{figure}
\def\w{4cm}
\def\dir{Figures/Camera/}
\includegraphics[width=\w]{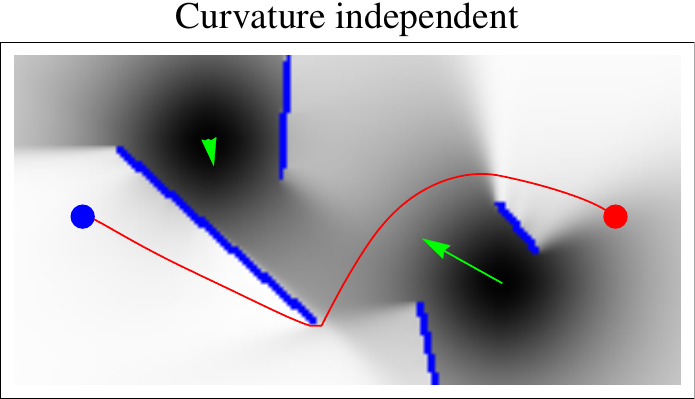}
\includegraphics[width=\w]{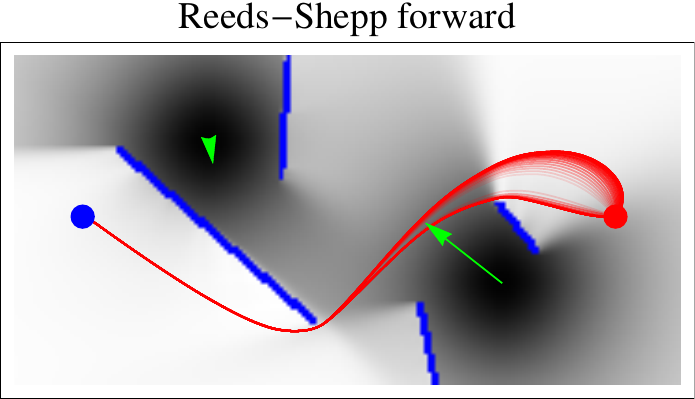}
\includegraphics[width=\w]{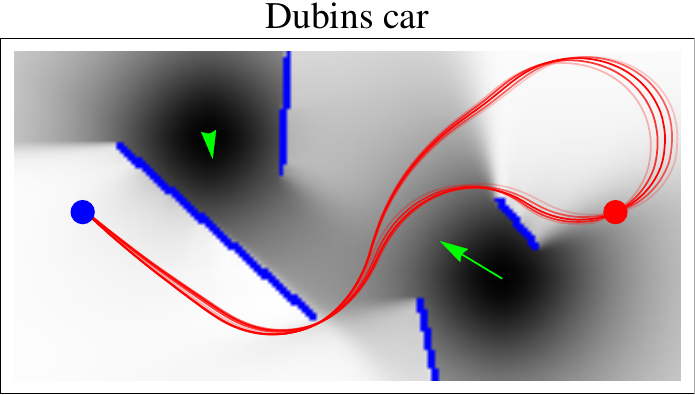}\\
\includegraphics[width=\w]{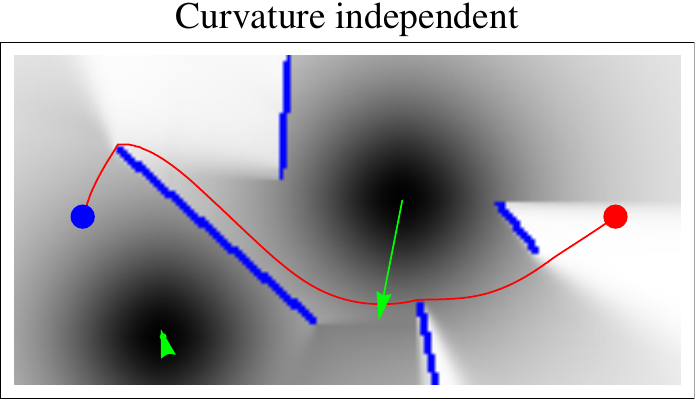}
\includegraphics[width=\w]{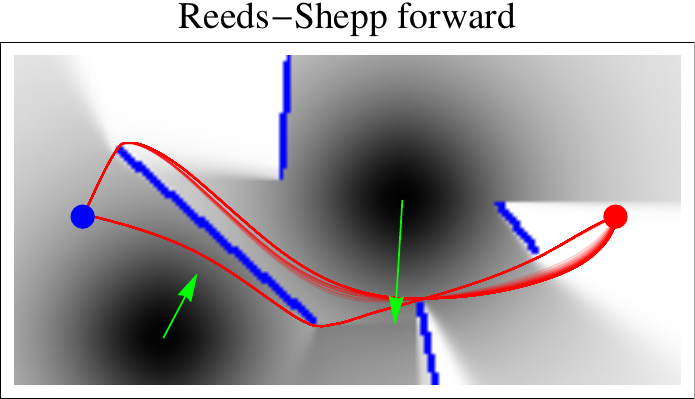}
\includegraphics[width=\w]{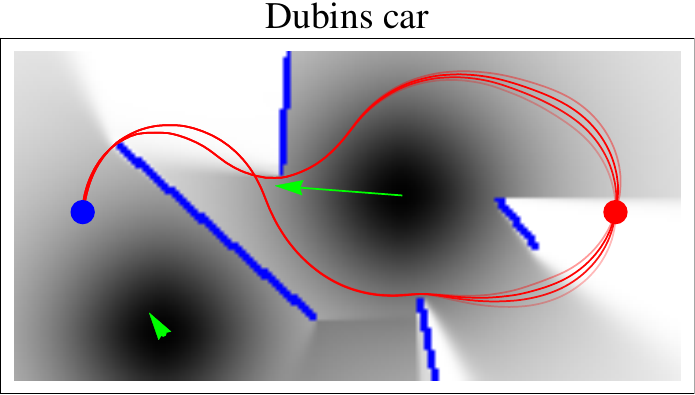}
\caption{Visual detection. Field of view of the cameras (black gradients), optimal furtive paths (red lines), local direction of improvement of the camera position (green arrows).}
\label{fig:Camera}
\end{figure}

\paragraph{Visual detection.} The first player places some cameras, which efficiency at detecting the second player decreases with distance and is blocked by obstacles, see Figure \ref{fig:Camera}. The cost function is 
\begin{equation}
    \cC_\xi(p,\dot p, \ddot p) = \cC_*(\kappa) \sum_{\substack{q \in \xi\\ [p,q] \subset R}} \frac 1 {\|q-p\|^2},
\end{equation}
where $\xi \in \Xi$ is a subset of $R$ with prescribed cardinality,  two in out experiments.
The green arrows on Fig \ref{fig:Camera} originate from the current (non optimal) camera position, and point in the direction of greatest growth $\nabla \gC(\xi,\Gamma)$ for the first player objective function.

\begin{figure}
\def\w{4cm}
\def\dir{Figures/Radar/}
\includegraphics[width=\w]{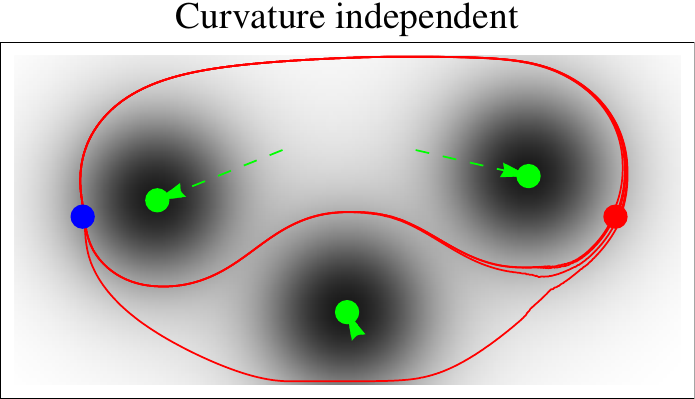}
\includegraphics[width=\w]{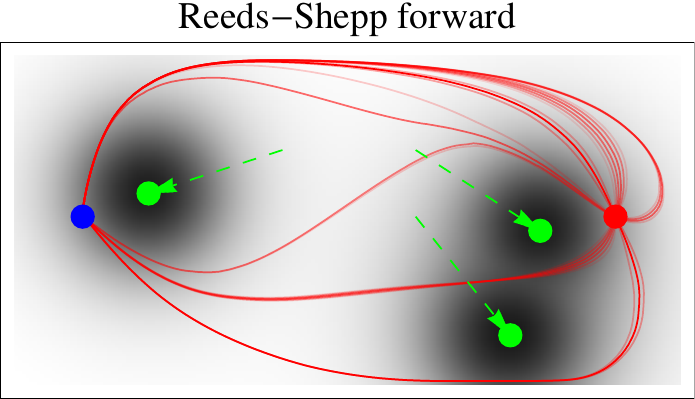}
\includegraphics[width=\w]{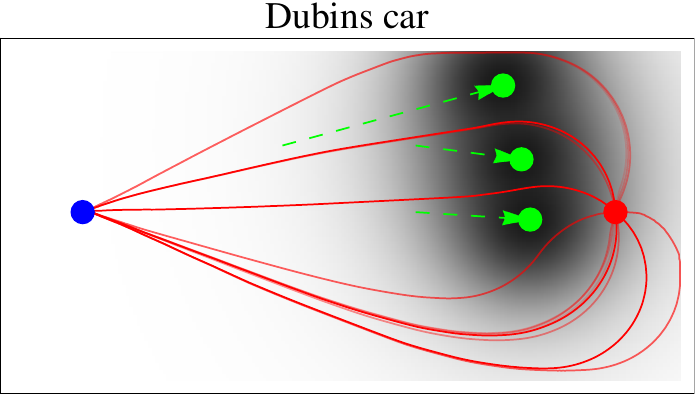}\\
\includegraphics[width=\w]{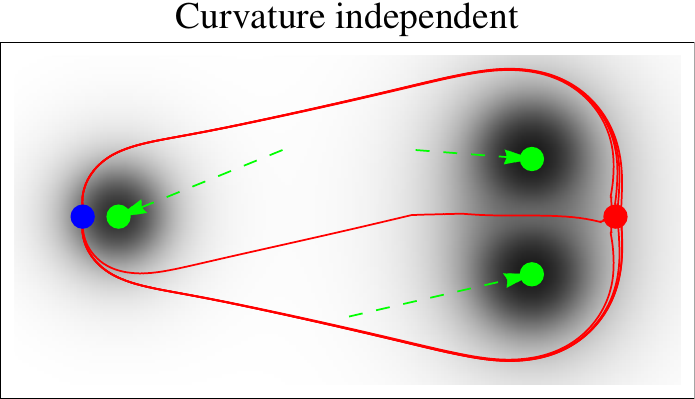}
\includegraphics[width=\w]{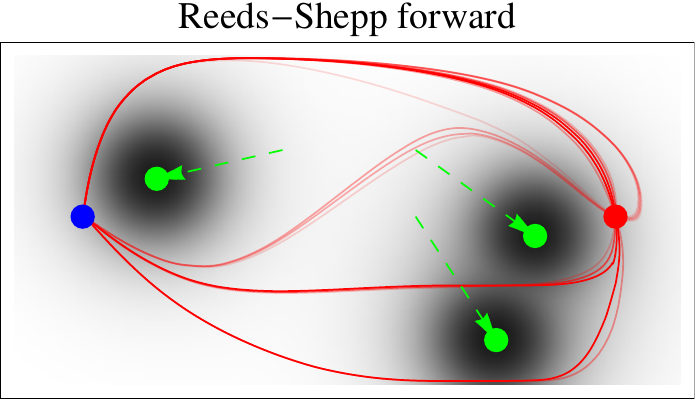}
\includegraphics[width=\w]{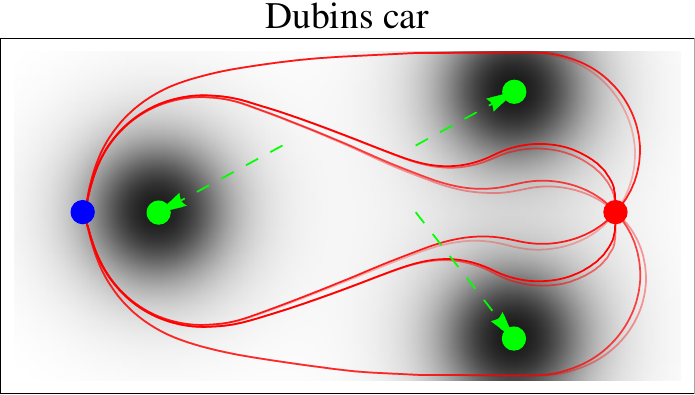}
\caption{Radar based detection. Optimal radar placement with an isotropic (top) or anisotropic (bottom) radar cross section.
}
\label{fig:Radar}
\end{figure}

\paragraph{Radar based detection.} The first player places some radars on the domain $R=[0,2]\times [0,1]$, here devoid of obstacles, and the second player has to fly by undetected.
The cost function is 
\begin{equation}
    \cC_\xi(p,\dot p, \ddot p) = \cC_*(|\ddot p|) \sqrt{\sum_{q \in \xi} \frac {\<\dot p,n_{\vec{pq}}\>^2+\delta^2\<\dot p,n^\perp_{\vec{pq}}\>^2}{\|p-q\|^4}}
\end{equation}
where $n_{\vec{pq}} := (q-p)/\|q-p\|$, and $\xi$ is a three element subset of $[0.4,1.6]$. The parameter $\delta$ is set to $1$ for an isotropic radar cross section (RCS), or to $0.2$ for an anisotropic RCS. In the latter where a plane showing its side to radar is five times less likely to be detected than a plane showing its nose or back, at the same position. Green arrows on Figure \ref{fig:Radar} point from the original position to the (locally) optimized position for player $\Xi$. At the equilibrium, several paths are optimal for player $\Gamma$, shown in red on Fig \ref{fig:Radar}.

\paragraph{Computational cost}
On a standard Laptop computer (2.7Ghz, 16GB ram), optimizing the second player objective, by solving a generalized eikonal equation, takes $\approx 1$s in the curvature dependent case, and $\approx 60$ times less in the curvature independent case thanks to the absence of angular discretization of the domain. Optimizing the first player objective takes $\approx 100$ $L-BFGS$ iterations, each one taking at most $8$s. For the stability of the minimization procedure, the problems considered were slightly regularized by the use of soft-minimum functions and by ``blurring'' the target keypoint over the $3\times 3$ box of adjacent pixels.

\section{Conclusion}

We have modeled a motion planning problem that take into account navigation constraints cand minimizing an anisotropic probability of detection, while computing the gradient of the value function related to the sensors location.
This model is thus useful for surveillance applications modeled as a two-player zero-sum game involving a target that tries to avoid detection.

%
% ---- Bibliography ----
%


\begin{thebibliography}{2}
%

\bibitem{BarbarescoMonnier:threatening}
F. Barbaresco and B. Monnier, {\it Minimal geodesics bundles by active contours: Radar application for computation of most threathening trajectories areas \& corridors}, 10th European Signal Processing Conference, Tampere, 2000, pp. 1--4.

\bibitem{Barbaresco:threatAniso}
F. Barbaresco, {\it Computation of most threatening radar trajectories areas and corridors based on fast-marching \& Level Sets}, IEEE Symposium on Computational Intelligence for Security and Defense Applications (CISDA), Paris, 2011, pp. 51--58.

\bibitem{Bardi1997:OptimalControl}
M. Bardi and I. Capuzzo-Dolcetta, {\it Optimal control and viscosity solutions of Hamilton-Jacobi-Bellman equations}, Bikhauser, 1997.

\bibitem{Benmansour2010:Derivatives}
F. Benmansour, G. Carlier, G. Peyré, and F. Santambrogio, {\it Derivatives with respect to metrics and applications: subgradient marching algorithm}, Numerische Mathematik, vol. 116, no. 3, pp. 357–381, May 2010.

\bibitem{Dubins1957:Car}
L. E. Dubins, {\it On curves of minimal length with a constraint on average curvature, and with prescribed initial and terminal positions and tangents}, Amer. J. Math., vol. 79, pp. 497–516, 1957.

\bibitem{Duits2017:ReedsShepp}
R. Duits, S.P.L. Meesters, J.-M. Mirebeau, J. M. Portegies,
{\it Optimal Paths for Variants of the 2D and 3D Reeds-Shepp Car with Applications in Image Analysis}, Preprint available on arXiv

\bibitem{Mirebeau2013:Diffusion}
J. Fehrenbach and J.-M. Mirebeau, {\it Sparse Non-negative Stencils for Anisotropic Diffusion,} Journal of Mathematical Imaging and Vision, pp. 1–25, 2013.

\bibitem{Mirebeau2014:FastMarching}
J.-M. Mirebeau, {\it Anisotropic Fast-Marching on cartesian grids using Lattice Basis Reduction}, SIAM J. Numer. Anal., vol. 52, no. 4, pp. 1573–1599, Jan. 2014.

\bibitem{Mirebeau2014:Minimal}
J.-M. Mirebeau, {\it Minimal stencils for discretizations of anisotropic PDEs preserving causality or the maximum principle}, SIAM J. Numer. Anal., vol. 54, no. 3, pp. 1582–1611, 2016.

\bibitem{Mirebeau2017:VR1}
J.-M. Mirebeau, {\it Anisotropic fast-marching on cartesian grids using Voronoi’s first reduction of quadratic forms}, in preparation

\bibitem{Mirebeau2017:Curv}
J.-M. Mirebeau, {\it Fast Marching methods for Curvature Penalized Shortest Paths}, in preparation

\bibitem{Mumford1994:Elastica}
D. Mumford, {\it Elastica and Computer Vision}, no. 31, New York, NY: Springer New York, pp. 491–506, 1994.

\bibitem{Rouy1992:Viscosity}
E. Rouy and A. Tourin, {\it A Viscosity Solutions Approach to Shape-From-Shading}, SIAM J. Numer. Anal., vol. 29, no. 3, pp. 867–884, Jul. 1992.

\bibitem{Schurmann2009:Quadratic}
A. Sch\"urmann, {\it Computational geometry of positive definite quadratic forms}, University Lecture Series, 2009.

\bibitem{Strode:multistatic}
C. Strode, {\it Optimising multistatic sensor locations using path planning and game theory}, IEEE Symposium on Computational Intelligence for Security and Defense Applications (CISDA), Paris, 2011, pp. 9--16.

\bibitem{Griewank2008:Derivatives}
A. Griewank, A. Walther,  {\it Evaluating derivatives: principles and techniques of algorithmic differentiation}. Society for Industrial and Applied Mathematics, 2008.

\end{thebibliography}
\end{document}